\input amstex
\documentstyle{amsppt}

\topmatter
\title 
Vanishing theorems and universal coverings of 
projective varieties
\endtitle
\rightheadtext{Vanishing theorems}
\author 
Fedor Bogomolov
\endauthor
\thanks
Partially supported by NSF grant DMS-9801591.
\endthanks
\affil 
Courant Institute for Mathematical Sciences 
\endaffil
\address
Courant Institute for Mathematical Sciences,
New York University
\endaddress
\email
bogomolo{\@}CIMS.NYU.EDU
\endemail
\abstract
This article contains a new argument which proves  vanishing of the
first cohomology for negative vector bundles over 
a complex projective variety if the rank of the bundle is
smaller than the dimension of the base.
Similar argument is applied to the 
construction of holomorphic functions on the universal covering of the complex
projective variety .
\endabstract
\endtopmatter

\document

\head 1. Introduction
\endhead

There are many flavors of  vanishing theorems for negative or semi-negative 
line bundles (see \cite{2} for a survey of the existing results). In contrast 
the vanishing of the cohomology of negative vector bundles is not 
so well understood. In this paper I propose a new approach leading to
a simple proof of the vanishing of the first cohomology of negative 
vector bundles whose rank is smaller than the dimension of the base.
This approach also provides some insight into the problem of quasi convexity
of the universal coverings of projective varieties.

For the results on Stein spaces used throughout the paper the reader
may wish to consult A.L.Onishchik's survey \cite{5}. Similarly to
\cite{5} the present  article emphasizes mainly the 
holomorphic and cohomological features of 
Stein spaces and keeps the references to differential geometry 
and plurisubharmonic
functions to a minimum. It is clear that most of the results
here can be strengthened by applying differential geometric tools
or by using some of the more elaborate properties of the vector bundles in
question.

The main idea of this article stems from Lemma~2.3 in 
my joint work \cite{1} with L.Katzarkov on the
fundamental groups of projective and symplectic manifolds.

\medskip

I am grateful to L.Katzarkov and T.Pantev for useful comments.

\head 2. A vanishing theorem for negative vector bundles
\endhead

Let $X$ be a normal irreducible complex projective variety. 

For a vector bundle $E \to X$ denote by $t(E) :=
\operatorname{Spec}(S^{\bullet}E^{\vee}) \to X$ its total space
scheme. Similarly let $\pi : {\Bbb P}(E) := {\bold
P}{\bold r}{\bold o}{\bold j}(S^{\bullet}E^{\vee}) \to X$ denote the
projectivization of $E$ and let ${\Cal O}_{{\Bbb P}(E)}(-1)
\subset \pi^{*}E$ denote the tautological line bundle. 

Recall (see e.g. \cite{3}) that a vector bundle $F$ is called {\it negative} if
the regular functions on $t(F)$ separate points and tangent directions
outside of the zero section. In other words $F$ is negative if the
natural affinization map
$$
a : t(F) \to \operatorname{Spec}({\Bbb C}[t(F)]) = 
\operatorname{Spec}(H^{0}(X,S^{\bullet}F^{\vee})) 
$$
contracts precisely the zero section.

\medskip

Our main result is the following theorem

\proclaim{Theorem 1} If $F$ is a negative vector bundle on $X$ with
$\operatorname{rk} F < \dim X$, then $H^1(X,F) = 0$.
\endproclaim

\demo{Proof} Put $d := \dim X$ and $r := \operatorname{rk} F$.
Assume that $H^1(X,F) \neq 0$. Then there is a cohomology class
$s\in H^1(X,F)$, $s\neq 0$. Let $\widetilde{F}$ be the corresponding 
extension of $F$ by ${\Cal O}_{X}$.
For a nonzero constant section $m\in {\Cal O}_{X}$ the preimage of $m$ in 
$\widetilde{F}$ is an affine bundle over $X$ modeled on $F$. We will 
denote it by $F_{s}$ since it does not depend on the choice of $m$.
Consider the closure $\overline{t(F)}$ of $t(F)$ and $\overline{F_s}$ of 
$F_s$ respectively inside the projective fibration
${\Bbb P}(\widetilde{F}\oplus {\Cal O}_{X})$. Note that both divisors 
at infinity $\overline{t(F)}\setminus t(F)$ and $\overline{F_s}\setminus F_s$
are isomorphic to ${\Bbb P}(F)$ with the same positive normal bundle $N \cong 
{\Cal O}_{{\Bbb P}(F)}(1)$.

By Grauert's criterion for ampleness \cite{3} the line bundles 
${\Cal O}_{F_{s}}({\Bbb P}(F))$ and \linebreak ${\Cal 
O}_{t(F)}({\Bbb P}(F))$  are both ample. Thus the linear system 
${\Cal O}_{F_{s}}({\Bbb P}(F))$ defines a
a projection $r : F_s \to \operatorname{Spec}({\Bbb C}[F_s]) =:  
F_s^r$ which is a proper morphism. Since by construction $r$ is birational
there are at most finitely many points $q_i\in F_s^r$ which have 
as preimages some positive dimensional compact subvarieties $X_i \subset 
F_s$.

Since each $X_i$'s is proper  the natural map $p_{|X_{i}} : X_i\to X$, 
where $p$ is the affine projection $p : F_s\to X$, will be finite 
on its image. In the next lemma we prove that $\dim_{{\Bbb C}} X_i <
\dim_{{\Bbb C}}  X = d$ for any $q_i$, which leads to a contradiction. 
Indeed  
$F_{s}^{r} = F_{s}/(\coprod_{i} X_{i})$ as a topological space (in the
classical topology) and so for
the reduced singular homology of $F_{s}$ with say ${\Bbb Q}$-coefficients we
have $\widetilde{H}_{i}(F_{s}^{r}, {\Bbb Q}) =
H_{i}((F_{s},\coprod_{i} X_{i}), {\Bbb Q})$. Now the long exact homology
sequence of the pair $(F_{s},\coprod_{i} X_{i})$ together with the
fact that $\coprod_{i} X_{i}$ is compact of complex dimension strictly
less than $d$ gives that $H_{2d}(F_{s},{\Bbb Q}) \cong
\widetilde{H}_{2d}(F_{s}^{r},  {\Bbb Q}) = H_{2d}(F_{s}^{r}, {\Bbb
Q})$. On the other hand since $F_{s}^{r}$ is 
an affine variety we know
that $F_{s}^{r}$ has the homotopy type of a simplicial complex 
of real dimension which does not exceed $\dim_{{\Bbb C}}F_{s}^{r} = 
\dim_{{\Bbb C}} t(F) = \dim_{{\Bbb C}} X + rk F = d + r$. This yields 
$H_{2d}(F_{s}^{r},{\Bbb Q}) = 0$ due to the
hypothesis $r < d$. This however is impossible since
$F_{s}$ is an affine bundle over $X$ and hence homotopically
equivalent to $X$ which in turn implies that $H_{2d}(F_{s},{\Bbb Q})$ 
contains a non-zero homology class - namely the
fundamental class of $X$. This gives us the desired contradiction and
finishes the proof of the theorem modulo the fact that $\dim_{{\Bbb
C}} X_{i} < \dim_{{\Bbb C}} X$ for all $i$.

We will derive this fact from the following general lemma.

\proclaim{Lemma 1} Let $X$ and $F$ be as in the statement of Theorem~1 and
let $y : Y \to X$ be a finite morphism. Then for any $0 \neq s \in
H^{1}(X,F)$ we have  $y^*s \neq 0$.
\endproclaim

\demo{Proof} Assume first that $X$ is a smooth compact curve. By
taking the Galois closure of $y : Y \to X$ and normalizing one gets 
a smooth finite Galois covering $y_G : X_G \to X$ with a finite
Galois group $G$ which factors a $y_G = y\circ y'$ with $y' : X_G \to Y$.
By assumption $y_G ^*s = 0$. Thus $y_G ^* \widetilde{F} \cong F\oplus
{\Cal O}_{X}$. Moreover since the property of being negative
is preserved by finite base change Proposition~4.3, \cite{3} we have
that $y_{G}^{*}F$ is a negative vector bundle on $X_{G}$. Consequently
$H^{0}(X_{G},y_{G}^{*}F) = 0$ and so
$H^{0}(X_{G},y_{G}^{*}\widetilde{F})$ injects into $H^{0}(X_{G},{\Cal
O}_{X_{G}}) \cong {\Bbb C}$. Therefore
$$
H^{0}(X_{G},y_{G}^{*}\widetilde{F}) = H^{0}(X_{G},{\Cal
O}_{X_{G}}) \cong {\Bbb C}
$$ 
and so $G$ acts trivially on the one dimensional space 
$H^{0}(X_{G},y_{G}^{*}\widetilde{F})$.  This implies that
$(y_{G*}y_{G}^{*}\widetilde{F})^{G} = (\widetilde{F}\otimes
y_{G*}{\Cal O}_{X_{G}})^{G} = \widetilde{F}$ has a nowhere vanishing
section. Hence the exact sequence $0 \to F \to \widetilde{F} \to {\Cal
O}_{X} \to 0$ splits and so $s= 0$ - a contradiction.

To finish the proof of the lemma it remains just to observe that
for any $X$ and any vector bundle $F$ we can find a sufficiently ample
divisor $D \subset X$ so that the cohomology group $H^{1}(X,F)$
injects into $H_{1}(D,F_{|D})$. Since negativity is preserved by
restrictions to subvarieties we can then find a smooth curve $C$
(cut out by finitely many general hyperplane sections) for which
the restriction of $s$ on $C$ is nonzero. 

The lemma is proven. \hfill \ \qed
\enddemo
 
We can now finish the proof of the theorem. Indeed since $X_{i}$ are
contracted by the affinization map $r$ we have that $p^{*}s_{|X_{i}} =
0$ for all $i$. Moreover since by construction $p_{|X_{i}} : X_{i} \to
X$ is finite onto its image Lemma~1 implies that $p_{|X_{i}}$ cannot
be surjective. Thus $\dim X_{i} < \dim X$ which concludes the proof of
the theorem. \hfill \ \qed
\enddemo

\remark{Remark 1} It is pretty easy to construct many negative bundles of 
rank $\geq \dim X$ with nontrivial first cohomology. One possible
construction goes as follows. Let $L$ be a very ample line bundle on
$X$ which gives an embedding $X \subset {\Bbb P}^{n}$. There 
is a surjective map $h: {\Cal O}_{X}^{\oplus n+1} \to L$ which 
defines a rank $n$ subbundle $\ker h = F \subset {\Cal O}_{X}^{\oplus
n+1}$. The vector bundle $F\otimes L^{-1} $ is a negative bundle with 
$H^1 (X,F\otimes L^{-1})\neq 0$. Namely there is a nontrivial element
$s \in H^1(X,F\otimes L^{-1})\neq 0$ which corresponds to the extension
$0 \to F\to {\Cal O}_{X}^{\oplus n+1} \to L \to 0$. To see that the
latter is non-split one only needs to observe that from the Euler
sequence we have an isomorphism $F\otimes L^{-1} \cong
\Omega^{1}_{{\Bbb P}^{n}|X}$ and so the extension class $s$ is just
the first Chern class of $L$.
\endremark

\remark{Remark 2} In the process of proving Theorem~1 we showed 
that for a negative vector bundle $F$ of arbitrary rank and any
element $s\in H^1(X,F)$ the union of all subvarieties $Y \subset X$
for which $s_{|Y} = 0$ is actually an algebraic subvariety of $X$.

This condition is somewhat reminiscent of  the Lang conjecture and
it seems likely there is some deeper relation between them. 
Because of this analogy we will formulate it in a separate lemma
\endremark

\proclaim{Lemma 2} Let $F$ be a negative bundle on a projective 
variety $X$ and let  $s\in H^1(X,F)$. Then there exists a finite set
of projective schemes  $\{X_i\}_{i \in I}$, $\dim X_i < \dim X$ 
and affine morphisms
$f_i : X_i \to X$ with the following property. 
Let $Y$ be an irreducible projective variety and let 
$f : Y \to X$ be a map with $f^*s = 0$ s.t. $f(Y)$ is not a point.
Then there exists a $i \in I$ $f$ factors as $f = f_i\circ g$ where  
$g : Y \to X_i$.
\endproclaim

\demo{Proof} Indeed in the notations of the proof of Theorem~1 we have
that $f^* \widetilde{F}$ splits into a direct sum $f^{*}F\oplus 
{\Cal O}_{Y}$ and therefore $f^*F_s$ is trivial as an affine bundle,
i.e. $F_{s}\times_{X} Y$ contains $Y$ as a closed subvariety. 
Consider the projection $p_{1} : F_{s}\times_{X} Y \to F_s$.
As we saw in the proof of Theorem~1 the affinization morphism of $F_s$ 
contracts only a finite set of proper subvarieties $X_{i} \in F_{s}$.
Hence $p_{1}(Y) = X_i$ for some $i$ and therefore $f = f_i\circ p_{1}$
which proves the lemma. \hfill \ \qed
\enddemo

\remark{Remark 3} As the proof shows the test variety $Y$ in the above lemma
can be any connected  variety without non constant holomorphic functions.
\endremark

\remark{Remark 4} If $E \to X$ is a positive vector bundle and $0 \neq
s \in H^{1}(X,E)$, then we 
can always find a variety $Y$ with
a surjective finite map $f : Y \to X$ for which $f^*s = 0$.
Indeed in this case the divisor ${\Bbb P}(E) \subset \overline{E}_{s}$
can be contracted and so by cutting down $\overline{E}_{s}/{\Bbb
P}(E)$ by a sequence of hyperplane sections we will eventually get a
subvariety $Y \subset E_{s}$ which has
dimension equal to $\dim X$. Since the natural projection $p : E_{s}
\to X$ is an affine morphism it follows that $p_{|Y} : Y \to X$ will
be a finite map and so we can take $f = p_{|Y}$.
\endremark

\head 3. Infinite covers of projective varieties
\endhead

Now we can apply the above results to the study of infinite covering
spaces of complex projective varieties.
 
\proclaim{Theorem 2} 
Let $F$ be a negative bundle over a complex projective variety $X$
and let $0 \neq s\in H^1(X,F)$ be a nontrivial cocycle. Let 
$Y$ be  a complex space with a map $f: Y\to X$  which is locally
finite and locally compact.
Assume that $f^*s = 0 $. Then there is a proper algebraic subset $Z\in X$
such that the holomorphic functions on $Y$ separate points on
$Y\setminus f^{-1}(Z)$.
\endproclaim 
 
\demo{Proof} 
Consider the variety $F_s$ as in the proof of Theorem~1. 
Holomorphic functions (even regular functions) separate points of
$F_s$ modulo a finite subset of algebraic subvarieties $X_i$ which
project finite to one  onto proper subvarieties in $X$.
Denote by $Z$ the union of the images of all $X_i$ in $X$.
Consider the manifold $F_{f^*s} = F_{s}\times_{X} Y$ which 
is an affine fibration  $p_Y : F_{f^*s}\to Y$ over $Y$.
  
The map $f$ induces a locally finite and locally compact map 
 $f' : F_{f^*s}\to F_s$. Thus the 
 holomorphic functions on $F_s$ locally separate points on $F_{f^*s}$ outside
of $f^{-1}(Z)$.
Since $f^*s = 0$ the corresponding extension of $F$ splits over $Y$.
Thus there is a section $r : Y \to F_{f^*s}$ of $f'$ and so the
 restriction of the holomorphic 
functions $F_{f^{*}s}$ to $r(Y)$ separate points outside of $r(f^{-1}Z)$.
This finishes the proof of the lemma. \hfill \ \qed
\enddemo

There are some immediate corollaries of this theorem.

\proclaim{Corollary 1} Let $X$, $F$ and $s$ be as in Theorem~2. Let
$f : \widetilde{X} \to X$ be the universal 
cover of $X$ and assume that $f^*s = 0$. Then the holomorphic functions on
$\widetilde{X}$ separate points outside a preimage of an algebraic 
subset $Z\in X$.
\endproclaim
\demo{Proof} Clear. \hfill \ \qed
\enddemo

\proclaim{Corollary 2} Let $X$, $F$ and $s$ be as in Theorem~2. Assume
furthermore that $F_{s}$ is an affine variety. Let $f : Y
\to X$ be any infinite unramified covering s.t. $f^{*}s = 0$. Then
$Y$ is Stein.
\endproclaim
\demo{Proof} Since any non-ramified covering of a Stein space is Stein
\cite{5} the assumption that $F_{s}$ is affine yields that
$F_{s}\times_{X} Y$ is Stein. On the other in the proof of Theorem~2
we saw that $Y \subset F_{s}\times_{X} Y$ is a closed analytic subset
and so $Y$ is Stein. \hfill \ \qed
\enddemo

\remark{Remark 5} The condition that $F_{s}$ is an affine variety can
be easily fulfilled in examples. For instance if $X$, $F$ and $s$ are
as in Remark~1 the variety $F_{s}$ is affine since by construction it
is a closed subvariety in the affine variety 
$${\Bbb P}^{n}\times {\Bbb
P}^{n \vee}\setminus \{(x,h) \in {\Bbb P}^{n}\times {\Bbb
P}^{n \vee}| x \in h \}.
$$
\
\endremark

Corollary~2 suggests that the result of Theorem~2 may also be
applicable to orbicoverings of $X$.

Let us first describe precisely the notion of orbicovering in the case of
a complex variety.  
 
Let $X$ be a complex variety and $S \subset X$ be a proper analytic
subset. Consider for any point $q\in S$ the local fundamental group 
$\pi_q = \pi_1(U(q)\setminus S)$ where $U(q)$ is a small 
ball in $X$ centered at $q$. Let $L \subset \pi_1(X\setminus S)$ 
be a subgroup with the property that $L\cap \pi_q$ is of finite index
in $\pi_{q}$ for all $q \in S$. Then the nonramified covering
of $ X\setminus S$ corresponding to $L$ can be naturally completed into
a normal complex variety $Y_{L}$ with a locally finite and locally compact
surjective map $f_L : Y_{L} \to X$. The map $f_L : Y_{L} \to X$ is
called an {\it orbicovering} of $X$ with a ramification set $S$. Now
we have the following

\proclaim{Corollary 3} Let $X$, $F$ and $s$ be as in Theorem~2. Assume
furthermore that $F_{s}$ is an affine variety. Let $f : Y
\to X$ be any orbicovering s.t. $f^{*}s = 0$. Then
$Y$ is Stein.
\endproclaim
\demo{Proof} Since every orbicovering of a Stein space is also Stein
(see Theorem~4.6 of \cite{5}) the proof is exactly the same as the
proof of Corollary~2. \hfill \ \qed
\enddemo

\remark{Remark 6} The prototype of Corollary~3 is Lemma~2.3 in our
joint work \cite{1}  with L.Katzarkov. 
\endremark

\bigskip

\noindent
In view of the previous results I would like to formulate a conjecture
which, I believe, should be the correct substitute  of the Shafarevich
holomorphic convexity conjecture. First we need the following definition.

\definition{Definition} Let $X$ be a smooth projective variety and let
$f : Y \to X$ be a nonramified covering of $X$. We will say that the
covering satisfies the property (N) (after M.Nori, T.Napier,
R.Narasimhan) if there exists a proper algebraic
subvariety $Z \subset X$ such that there exists normal complex space
$M$ satisfying

\smallskip

(a) \quad The holomorphic functions on $M$ separate points.

\smallskip

(b) \quad There exists a proper map with connected fibers 
$g :Y\setminus f^{-1}(Z) \to M$

(c) \quad Holomorphic functions on $Y$ separate the fibers
of $g$ after restriction to $Y\setminus f^{-1}(Z) $
 
\enddefinition

The conjecture now reads

\proclaim{Conjecture 1} Let $X$ be a smooth 
compact projective variety. Then he universal cover of $X$  satisfies 
the property (N).
\endproclaim

\remark{Remark 7} We have shown above that if there exists a negative
vector bundle $F \to X$ with a cohomology class $s \in H^1(X,F)$ which
becomes trivial on $\widetilde{X}$, then the conjecture is true and
moreover $\dim M = \dim X$. Thus the identification on $\widetilde{X}$
of the two "infinitesimally" close bundles $F\oplus {\Cal O}$ and
$\widetilde{F}$ implies the property (N) for the universal covering.

It is worth to point out that in many of the cases for which the
Shafarevich conjecture is known the proof relies on the comparison of
two vector bundles on $X$ which become equal when they are pulled back
to $\widetilde{X}$ \cite{3}. For example the theorem of M.Gromov
 uses the fact that some positive line
bundle $E \to X$ becomes trivial when pulled back to the universal
cover. The theorem of L.Katzarkov  \cite{4} establishes the 
holomorphic convexity of
$\widetilde{X}$  for projective surface $X$ under
the assumption of the existence of an almost faithful 
linear representation of $\pi_{1}(X)$. In this  case 
all the bundles  on $X$ corresponding to the representations 
of the fundamental group
 of the same dimension are becoming equal on $\widetilde{X}$.
\endremark

\remark{Remark 8} In the case of surfaces the space 
$\widetilde{X}$ can be obtained as a union of two rather simple 
Stein manifolds with Stein intersection. This implies that the structure 
of the space of the moduli space of vector bundles on $\widetilde{X}$ 
in this case is somewhat similar to the structure of the moduli space of 
vector bundles on a  curve. Namely any bundle of rank greater than
 $2$ has a complete flag 
of subbundles, thus reducing the K-group $K_{0}(\widetilde{X})$ to
$\operatorname{Pic}(X)\times {\Bbb Z}$.
In particular one expects that many different bundles on $X$ coincide 
after lifting to $\widetilde{X}$. 
\endremark

\enddocument